\documentclass{amsart}

\usepackage[T1]{fontenc}

\usepackage[ngerman,turkish,american]{babel}
\usepackage{a4wide}
\usepackage{amssymb}
\usepackage[latin1]{inputenc}
\usepackage[T1]{fontenc}
\usepackage{graphicx}
\usepackage{color}
\usepackage[backref,hyperref]{hyperref}

\renewcommand{\phi}{\varphi}
\newcommand{\C}{{\mathbb{C}}}
\newcommand{\R}{{\mathbb{R}}}
\newcommand{\Z}{{\mathbb{Z}}}
\newcommand{\N}{{\mathbb{N}}}

\renewcommand{\epsilon}{\varepsilon}
\renewcommand{\theta}{\vartheta}
\renewcommand{\S}{{\mathbb{S}}}
\newcommand{\D}{{\mathbb{D}}}

\newcommand{\Reeb}{{X_\mathrm{Reeb}}}

\newcommand{\z}{{\mathbf{z}}}
\newcommand{\lie}[1]{{\mathcal{L}_{#1}}}

\newcommand{\abs}[1]{{\lvert #1\rvert}}

\newcommand{\ba}{\begin{array}}
\newcommand{\ea}{\end{array}}

\DeclareMathOperator{\lcm}{lcm}

\DeclareMathOperator{\Vol}{Vol}

\theoremstyle{plain}
\newtheorem{theorem}{Theorem}
\newtheorem{lemma}[theorem]{Lemma}

\theoremstyle{remark}
\newtheorem{remark}{Remark}

\theoremstyle{definition}
\newtheorem*{defi}{Definition}

\begin{document}
\bibliographystyle{amsalpha}

\title{Brieskorn manifolds as contact branched covers of spheres}

\selectlanguage{turkish}
\author{Fer\.{\i}t Öztürk}
\address{Bo\u{g}az\.{\i}ç\.{\i} Ün\.{\i}vers\.{\i}tes\.{\i}, Department of Mathematics, TR-34342
  Bebek, \.Istanbul, Turkey}
\selectlanguage{american}
\email{ferit.ozturk@boun.edu.tr}

\selectlanguage{ngerman}
\author{Klaus Niederkrüger}
\address{Mathematisches Institut, Universität zu Köln\\Weyertal
  86-90\\D-50.931 Köln\\Federal Republic of Germany}
\selectlanguage{american}
\email{kniederk@math.uni-koeln.de}

\begin{abstract}
  We show that Brieskorn manifolds with their standard contact
  structures are contact branched coverings of spheres. This covering
  maps a contact open book decomposition of the Brieskorn manifold
  onto a Milnor open book of the sphere.

\end{abstract}

\maketitle

\section{Introduction}
Brieskorn manifolds have been an interesting source of examples. In the
field of topology many exotic spheres were found to be such manifolds, but
also in contact geometry they have provided a rich family of examples,
especially in dimensions larger than~$3$.

It has been known for a long time that a Brieskorn manifold
$\Sigma(a_0,\dotsc,a_n)\subset\C^{n+1}$ is an $a_0$-fold cyclic covering of
the unit sphere $\S^{2n-1}\subset \C^n$ branched along the
$(2n-3)$-dimensional Brieskorn manifold $\Sigma(a_1,\dotsc,a_n)$.  In this
article, we show that this is not only true as smooth manifolds but also in
the contact category.

Furthermore, there is for every $(2n-3)$-dimensional Brieskorn
manifold $B:=\Sigma(a_1,\dotsc,a_n)$ a (so-called) Milnor open book on
$\S^{2n-1}$ that has $B$ as its binding.  This open book decomposition
can be pulled back by the cyclic branched covering to the Brieskorn manifold
$\Sigma(a_0,\dotsc,a_n)$. 
In this way it is possible to show that the
open book of a Brieskorn manifold can be described in an abstract way
by using a Milnor open book and taking a  power of its
monodromy map. One of the goals of this article is to show that 
the canonical contact structure on a Brieskorn
manifold is supported by that open book (see Section \ref{ob and ctct} 
for the definition).

Below we will first state what a contact branched covering is.
Lemma~\ref{branched cover gives open book} then shows that the contact
structure of the contact branched covering is supported by the open book
induced by the compatible open book in the base. Theorem~\ref{milob} proves
that there is a contact structure isotopic to the standard one on
$\S^{2n-1}$ which is supported by the Milnor open book of $\S^{2n-1}$ with
binding a Brieskorn manifold $\Sigma(a_1,\dotsc,a_n)$. Finally,
Theorem~\ref{brieskorn is ctct branched cover} shows that the contact
structure of the $a_0$-fold contact branched covering of $\S^{2n-1}$ as in
Lemma~\ref{branched cover gives open book} is isotopic to the standard
contact structure on the Brieskorn manifold $\Sigma(a_0,a_1,\dotsc,a_n)$.

\section{Contact branched coverings}

Branched coverings for contact $3$-manifolds were first considered by
Gonzalo in \cite{Gonzalo}. He used them to reprove the existence of a
contact structure on any oriented $3$-manifold. His methods used local
charts and were adapted to his special situation. Geiges showed later that
a branched covering of a contact manifold of any dimension admits under
very natural conditions a contact structure \cite{Geiges_Construction}.
Below we will give a definition of contact branched covers, which coincides
essentially with Geiges' construction, and show that up to isotopy it is
independent of any choices.

Let $(N,\alpha)$ be a contact manifold, and let $f:\,M\to N$ be a
branched covering.  The pull-back form $f^*\alpha$ fails to be contact 
on $M$, because by definition $\dim(\ker df) = 2$ along the branching locus.  This
problem can be fixed though by perturbing $f^*\alpha$ slightly.

\begin{lemma}\label{contact branched covering}
  Let $f:\,M\to N$ be a covering branched along $B\subset N$ such that
  $(N,\alpha)$ and $(B,\left.\alpha\right|_{TB})$ are contact manifolds.
  There exists a $1$-form $\gamma$ on $M$ with $\left.d\gamma\right|_{\ker
    df} > 0$ along $B$ ($\ker df$ is naturally oriented, because
  $f^*\alpha$ gives an orientation for $M$ and $f^{-1}(B)$) such that
  $$
  f^*\alpha + \epsilon\,\gamma
  $$
  is a contact form on $M$ for any sufficiently small~$\epsilon>0$.
  
  Any contact form $\beta_1$ on $M$ is isotopic to $f^*\alpha +
  \epsilon\,\gamma$ if it lies in a smooth family of $1$-forms $\beta_t$
  with $t\in[0,1]$ such that $\beta_0 = f^*\alpha$, and for which $\beta_t$
  is contact for all $0<t\le 1$, and for which
  $\left.d\widetilde\gamma\right|_{\ker df} > 0$, where we have
  set~$\widetilde\gamma = \left.\dot\beta_t\right|_{t=0}$.
\end{lemma}

\begin{defi}
  $f:\,M\to N$ together with the contact structure given above is
  called the \textbf{contact branched covering} of $(N,\alpha)$ along
  $(B,\left.\alpha\right|_{TB})$.
\end{defi}

\begin{proof}
  The existence of such a form $\gamma$ was proved in \cite{Geiges_Construction},
  and the uniqueness of the contact structure can be shown in a similar
  way.  For completeness though, here is the argument: Consider the Taylor
  expansion of $\beta_t$ at $t=0$
  \begin{align*}
    \beta_t &= f^*\alpha + t\,\widetilde\gamma + \mathcal{O}(t^2)\;.
  \end{align*}  
  We will use this $1$-form at time $t_0=\epsilon>0$, where $\epsilon$ will
  be chosen below.  We can form the linear interpolation between
  $\beta_\epsilon$ and $f^*\alpha + \epsilon\,\gamma$ to define the family
  of $1$-forms
  \begin{align*}
    \alpha_s := f^*\alpha + \epsilon\,(s\,\gamma + (1-s)\,
    \widetilde\gamma) + (1-s)\,\mathcal{O}(\epsilon^2)\quad\text{ for
      $s\in[0,1]$.}
  \end{align*}
  The contact condition for this family becomes
  \begin{align*}
    \begin{split}
      \alpha_s\wedge \bigl(d\alpha_s\bigr)^n & := f^*\left(\alpha\wedge
        \bigl(d\alpha\bigr)^n\right) + \epsilon\,(s\,\gamma +
      (1-s)\, \widetilde\gamma)\wedge f^*\bigl(d\alpha\bigr)^n  \\
      & \qquad + \epsilon
      n\,f^*\bigl(\alpha\wedge\bigl(d\alpha\bigr)^{n-1}\bigr)\wedge
      (s\,d\gamma + (1-s)\, d\widetilde\gamma) + \mathcal{O}(\epsilon^2)\;.
    \end{split}
  \end{align*}
  
  On the branching locus, the first two terms vanish; the third one is
  positive for all $s\in[0,1]$ by our assumptions, and by choosing
  $\epsilon>0$ small enough it dominates the
  $\mathcal{O}(\epsilon^2)$-part.  By continuity there is an open
  neighborhood $U$ of $f^{-1}(B)$ where the sum of all terms containing an
  $\epsilon$-factor is positive for any sufficiently small $\epsilon > 0$.
  The pull-back $f^*\left(\alpha\wedge \bigl(d\alpha\bigr)^n\right)$ is
  positive on the compact set $M-U$, and is thus always larger than
  $C\,\Vol_M$ for some $C>0$. We can achieve that the $\epsilon$-terms (by
  choosing $\epsilon$ still smaller if necessary) are never smaller than
  $-C\,\Vol_M$. For any sufficiently small $\epsilon > 0$, it follows that
  $\alpha_s\wedge d\alpha_s^n > 0$, and thus the corresponding contact
  structures are isotopic by Gray stability.
\end{proof}

Note that the definition of a contact branched covering is in analogy
with the definition of a symplectic branched covering  \cite{Aur}.
Furthermore there is the concept of canonicity of the  structure in the  
symplectic setting too (see \cite[Proposition~3]{Aur}).

\section{Open books and contact structures}
\label{ob and ctct}

The following definitions are taken from \cite{Giroux}.

\begin{defi}
  An \textbf{open book} on a closed manifold $M$ is given by a
  codimension-$2$ submanifold $B\hookrightarrow M$ with trivial normal
  bundle, and a bundle $ \theta:\,(M-B)\to \S^1$.  The neighborhood of
  $B$ should have a trivialization $\D^2\times B$, where the angle
  coordinate on the disk agrees with the map~$\theta$.
  
  The manifold $B$ is called the \textbf{binding} of the open book and
  a fiber $P=\theta^ {-1}(\varphi_0)$ is called a \textbf{page}.
\end{defi}

\begin{remark}
  The open set $M-B$ is a bundle over $\S^1$, hence it is diffeomorphic to
  the mapping torus $P_\Phi := \R\times P/\sim$, where $\sim$ identifies
  $(t,p)\sim(t+1,\Phi(p))$ for some diffeomorphism $\Phi$ of $P$. Since the
  neighborhood of the binding has the standard form described above, we can
  assume that $\Phi$ is equal to the identity in some small neighborhood of
  the binding.  By glueing $\D^2\times B \cong \D^2 \times \partial P$ onto
  $P_\Phi$ in the obvious way, we obtain a manifold diffeomorphic to $M$.
\end{remark}

\begin{defi}\label{compatiblebook}
  A contact structure $\xi$ on $M$ is said to be \textbf{supported by
    an open book $(B,\theta)$} of $M$, if there is a contact form
  $\alpha$ with $\xi=\ker\alpha$ such that
  \begin{enumerate}
  \item $(B,\alpha|_{TB})$ is a contact manifold.
  \item For every $s\in \S^1$, the page $P:=\theta^{-1}(s)$ is a
    symplectic manifold with symplectic form $d\alpha$.
  \item Denote the closure of a page $P$ in $M$ by $\overline{P}$. The
    orientation of $B$ induced by its contact form $\alpha|_{TB}$
    should coincide with its orientation as the boundary of
    $(\overline{P},d\alpha)$.
  \end{enumerate}
  Such a contact form is said to be \textbf{adapted} to $(B,\theta)$.
\end{defi}

\begin{remark}\label{orientation condition for connected binding}
  Note that if the binding is connected, point~(3) of the definition
  above holds automatically, because
  \begin{align*}
    0 & < \int_P (d\alpha)^n = \int_B \alpha\wedge (d\alpha)^{n-1}\;,
  \end{align*}
  by Stokes theorem. Hence the orientation of $B$ as boundary of $P$
  agrees with the one given by the contact form.
\end{remark}

\begin{lemma}\label{branched cover gives open book}
  Let $(N,\alpha)$ be a contact manifold that has an open book
  decomposition $(B,\theta)$ supporting $\alpha$.  The $k$-fold cyclic
  covering $f:\,M\to N$ branched over $B$ exists, and is a contact manifold
  adapted to the open book decomposition $(f^{-1}(B),\sqrt[k]{\theta\circ f})$.
\end{lemma}
\begin{proof}
  Note that $N-B$ can be written by the remark above as $P_\Phi = \R\times
  P/\sim$, where $\sim$ identifies $(t,p)$ with $(t+1,\Phi(p))$ for some
  diffeomorphism $\Phi$ of the page $P$ that is the identity in a small
  neighborhood of $\partial P$.
  
  Construct $M$ as the mapping torus $P_{\Phi^k} = \R\times P/\sim_k$,
  where $\sim_k$ identifies $(t,p)$ with $(t+1,\Phi^k(p))$ for the
  diffeomorphism $\Phi$ on $P$. At the boundary the mapping torus is still
  diffeomorphic to $\S^1\times (-\epsilon,0]\times \partial P$ such that we
  can glue in $\D^2\times B$ to obtain a closed manifold $M$.

  Define the projection $f:\,M\to N$ of the branched covering piecewise:
  \begin{align*}
    \begin{array}{ccccc}
      M & \cong & P_{\Phi^k} & \cup_{\S^1\times\partial P} & \D^2\times B
      \\
      f\downarrow &  & \downarrow f_1 & & \downarrow f_2   \\
      N & \cong & P_\Phi & \cup_{\S^1\times\partial P} & \D^2\times B \\
    \end{array}
  \end{align*}
  The map $f_1:\,P_{\Phi^k} \to P_\Phi$ is given by $f_1([t,p]) = [kt,p]$,
  and the map $f_2:\, \D^2\times B \to \D^2\times B$ is given by
  $f_2(re^{i\phi},p) = (g(r) e^{ik\phi},p)$, where $g(r)$ is a smooth
  strictly increasing function on $\R_{\ge 0}$ that is equal to $r^k$ close
  to zero and equal to $r$ for $r>\delta$ with $\delta>0$ very small. Then
  it is clear that $f$ defines a branched covering.
  
  It is clear by Lemma~\ref{contact branched covering} that $M$ supports a
  contact structure compatible with $f$. The contact form on $M$ is
  obtained by taking the pull-back $f^*\alpha$ and adding a small $1$-form
  $\gamma$ such that $\left.d\gamma\right|_{\ker df} > 0$. This $\gamma$
  can be chosen to be of the form $\gamma = \epsilon r^2\rho(r)\,d\phi$ on
  $\D^2\times B$.
  
  It is also clear that $(f^{-1}(B),\sqrt[k]{\theta\circ f})$ is an open
  book decomposition of $M$. Since $d\gamma$ vanishes, when restricted to
  any page, it follows that $\alpha + \gamma$ is supported by this open
  book.
\end{proof}

\section{Brieskorn manifolds and their canonical contact
structures} \label{briesk}

Before talking about Brieskorn manifolds, we will briefly collect some
facts about the sphere: Assume $\S^{2n-1}$ to be embedded in the
standard way in $\C^n$.  We will denote the points of $\C^n$ by $\z =
(z_1,\dotsc,z_n)$.  The standard contact form on the sphere is
\begin{align*}
  \alpha_\mathrm{std} &= \frac{i}{2} \sum_{j=1}^n (z_j\,d\bar z_j- \bar z_j\,dz_j)\;.
\end{align*}

\begin{lemma}
  \label{isot}
  The $1$-form
  \begin{align*}
    \beta &= \frac{i}{2}\sum_{j=1}^n a_j(z_j\,d\bar z_j - \bar
    z_j\,dz_j)\;,
  \end{align*}
  with $a_j\in\N$, is isotopic to the standard contact form on
  $\S^{2n-1}\subset\C^n$.
\end{lemma}
\begin{proof}
  The proof works by taking the linear interpolation between $\beta$ and
  $\alpha_\mathrm{std}$, and checking that all forms in the family are
  contact. This allows us to use Gray stability.
\end{proof}

Now, we will explain what a Brieskorn manifold is.
Let $f:\,\C^{n+1}\to\C$ be a polynomial of the form
\begin{align*}
  f(z_0,z_1,\dotsc,z_n) & = z_0^{a_0} + \dotsb + z_n^{a_n}\;,
\end{align*}
with fixed numbers $a_0,\dotsc,a_n\in \N$.  It is easy to see that the
variety $V_f := f^{-1}(0)$ has a single isolated singularity at
$(0,\dotsc,0)$. Outside the origin, the equation describes a smooth
submanifold of codimension $2$, because the matrix
\begin{align*}
  \left(
    \begin{matrix}
      \partial f & \bar \partial f \\
      \partial \bar f & \bar \partial \bar f
    \end{matrix}
  \right) &=
  \left(\begin{matrix}
      a_0 z_0^{a_0-1} & \hdots & a_n z_n^{a_n-1} & 0 & \hdots & 0 \\
      0 & \hdots & 0  & a_0 \bar z_0^{a_0-1} & \hdots & a_n \bar z_n^{a_n-1}
  \end{matrix}\right)
\end{align*}
has full rank.

\begin{defi}
  The \textbf{Brieskorn manifold} $\Sigma(a_0,\dotsc,a_n)$ is defined
  as the intersection
  \begin{align*}
    \Sigma(a_0,\dotsc,a_n) := V_f \cap \S^{2n+1}\;.
  \end{align*}
\end{defi}

This set is, as its name suggests, a manifold.  This can be easily seen by
noting that $V_f$ is transverse to $\S^{2n+1}$.  Since the sphere has
codimension $1$, it is enough to find a vector field $Z$ on $V_f$, which is
everywhere transverse to the sphere.  The $\R$-action
\begin{align*}
  \R\times \C^{n+1} & \to \C^{n+1} \\
  (z_0,\dotsc,z_n) & \mapsto
  \left(e^{t/a_0}z_0,\dotsc,e^{t/a_n}z_n\right)
\end{align*}
restricts to the variety $V_f$, and its infinitesimal generator
\begin{align*}
  Z(z_0,\dotsc,z_n) &= \left(\frac{z_0}{a_0},\dotsc,\frac{z_n}{a_n}\right)
\end{align*}
is always transverse to the sphere, because
\begin{align*}
\lie{Z}(\abs{z_0}^2 + \dotsb + \abs{z_n}^2-1) &= \frac{1}{a_0}
  \abs{z_0}^2 + \dotsb + \frac{1}{a_n} \abs{z_n}^2 \ne 0\;.
\end{align*}

\begin{figure}[htbp]
\begin{picture}(0,0)%
\includegraphics{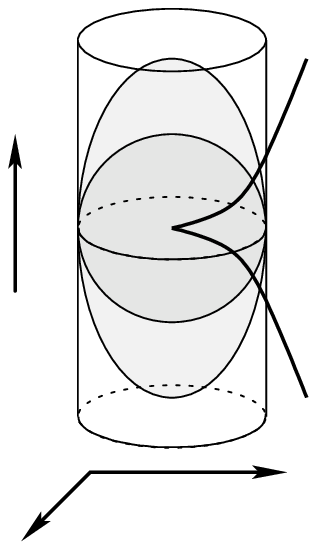}%
\end{picture}%
\setlength{\unitlength}{4144sp}%
\begingroup\makeatletter\ifx\SetFigFont\undefined%
\gdef\SetFigFont#1#2#3#4#5{%
  \reset@font\fontsize{#1}{#2pt}%
  \fontfamily{#3}\fontseries{#4}\fontshape{#5}%
  \selectfont}%
\fi\endgroup%
\begin{picture}(1639,2539)(72,-2063)
\put(1711,119){\makebox(0,0)[lb]{\smash{\SetFigFont{12}{14.4}{\familydefault}{\mddefault}{\updefault}{\color[rgb]{0,0,0}$V_f$}%
}}}
\put(406,344){\makebox(0,0)[lb]{\smash{\SetFigFont{12}{14.4}{\familydefault}{\mddefault}{\updefault}{\color[rgb]{0,0,0}$C_0$}%
}}}
\put(656,-1936){\makebox(0,0)[lb]{\smash{\SetFigFont{12}{14.4}{\familydefault}{\mddefault}{\updefault}{\color[rgb]{0,0,0}$(z_1,\ldots,z_n)$}%
}}}
\put( 72,-618){\makebox(0,0)[lb]{\smash{\SetFigFont{12}{14.4}{\familydefault}{\mddefault}{\updefault}{\color[rgb]{0,0,0}$z_0$}%
}}}
\end{picture}
  \caption{The manifolds $\Sigma(a_0,\dotsc,a_n)$ and
    $\widetilde\Sigma(a_0,\dotsc,a_n)$ are obtained by intersecting $V_f$
    with different hypersurfaces.}
  \label{bild von brieskorn und zylinder}
\end{figure}

In the rest of the article, we make extensive use of a related manifold:
Instead of taking the intersection between $V_f$ and a sphere, define
\begin{align*}
  \widetilde \Sigma (a_0,\dotsc,a_n) & := V_f \cap C_0\;,
\end{align*}
where $C_0$ is the spherical cylinder given by
\begin{align*}
  C_0 := \C\times \S^{2n-1} = \bigl\{(z_0,z_1,\dotsc,z_n)\bigm| (z_1,\dotsc,z_n)\in
  \S^{2n-1}\bigr\}\;.
\end{align*}
As above it is easy to check that this set is a manifold, because for
the defining equation of $C_0$, we obtain
\begin{align*}
  \lie{Z}(\abs{z_1}^2 + \dotsb + \abs{z_n}^2 - 1) &= \frac{1}{a_1}
  \abs{z_1}^2 + \dotsb + \frac{1}{a_n} \abs{z_n}^2 \ne 0\;.
\end{align*}

The Brieskorn manifold is of course diffeomorphic to $\widetilde
\Sigma(a_0,\dotsc,a_n)$ (see Figure~\ref{bild von brieskorn und
  zylinder}). In fact, let
\begin{align*}
  R_s &:= s\,\abs{z_0}^2 + \abs{z_1}^2 + \dotsb + \abs{z_n}^2\;,
\end{align*}
then we can define a family of submanifolds $\Sigma_s$ with
$s\in[0,1]$ by
\begin{align*}
  \Sigma_s := V_f \cap R_s^{-1}(1)\;,
\end{align*}
where $\Sigma_1$ is equal to $\Sigma(a_0,\dotsc,a_n)$ and
$\Sigma_0$ is equal to $\widetilde\Sigma(a_0,\dotsc,a_n)$.

\begin{lemma}
  There is an isotopy $\Phi_s$ in $V_f$ between
  $\Sigma(a_0,\dotsc,a_n)$ and $\Sigma_s$.
\end{lemma}
\begin{proof}
  Consider the $\R$-flow above, but let the time-parameter depend on
  the point that is being mapped, i.e.\ consider the map
  \begin{align*}
    \Phi_s:\,(z_0,\dotsc,z_n) & \mapsto
    \left(e^{T/a_0}z_0,\dotsc,e^{T/a_n}z_n\right)\;,
  \end{align*}
  where $T = T(z_0,\dotsc,z_n;s)$ is a function with the following
  properties: For a point $(z_0,\dotsc,z_n) \in
  \Sigma(a_0,\dotsc,a_n)$, we want its image to lie in $\Sigma_s$,
  hence the equation
  \begin{align*}
    1 & = s\,\abs{e^{T/a_0}z_0}^2 + \abs{e^{T/a_1}z_1}^2 + \dotsb +
    \abs{e^{T/a_n}z_n}^2 = s\,e^{2T/a_0}\abs{z_0}^2 +
    e^{2T/a_1}\abs{z_1}^2 + \dotsb + e^{2T/a_n}\abs{z_n}^2
  \end{align*}
  needs to hold. For any point $(z_0,\dotsc,z_n)$ there is a unique
  solution $T(z_0,\dotsc,z_n;s) \ge 0$, because the right hand side of the
  equation is a strictly increasing continuous function in $T$ that takes a
  value less than 1 for $T=0$.

  To prove that the map $\Phi_s$ is a bijection, construct a map
  $\widetilde\Phi_s$ analogously to the one above, which maps
  $\Sigma_s$ into $\Sigma(a_0,\dotsc,a_n)$. It is easy to see that
  these maps are mutually inverse.

  That $\Phi_s$ is smooth follows from the fact that $T$ is, and this is
  proved checking the inequality:
\begin{align*}
    \frac{d}{dT}\left(se^{2T/a_0}\abs{z_0}^2 + e^{2T/a_1}\abs{z_1}^2 + \dotsb +
      e^{2T/a_n}\abs{z_n}^2 - 1\right) > 0\;,
  \end{align*}
  which allows us to apply the implicit function theorem. The map $\Phi_s$
  is a bijective local diffeomorphism between closed manifolds, hence it is
  a diffeomorphism.
\end{proof}

\begin{lemma}
\label{alfa_s}
  For every $\Sigma_s$ with $s\in (0,1]$, the corresponding $1$-form
  \begin{align*}
    \alpha_s := \frac{i}{2} \, \Bigl( s a_0\,(z_0\,d\bar z_0 - \bar z_0\,dz_0) +
    a_1\,(z_1\,d\bar z_1 - \bar z_1\,dz_1) + \dotsb + a_n\,(z_n\,d\bar
    z_n - \bar z_n\,dz_n) \Bigr)
  \end{align*}
  is a contact form, and by Gray stability it follows that every $\Sigma_s$
  (with the exception of $\Sigma_0 = \widetilde\Sigma(a_0,\dotsc,a_n)$) is
  contactomorphic to $\Sigma(a_0,\dotsc,a_n)$.
\end{lemma}
\begin{proof}
  A long but trivial calculation yields
  \begin{align*}
    \begin{split}
      \alpha_s\wedge d\alpha_s^{n-1}\wedge dR_s\wedge df \wedge d\bar
      f & = \Bigl(s\bar f \sum_{j=0}^n a_jz_j^{a_j} + sf \sum_{j=0}^n
      a_j\bar z_j^{a_j}\\
      & \qquad -2a_0 R_s\abs{z_0}^{2(a_0-1)} -2s R_s \sum_{j=1}^n a_j
      \abs{z_j}^{2(a_j-1)} \Bigr)\, \Omega \;,
  \end{split}
\intertext{with $\Omega := i^n/2(n-1)!\, a_0\cdots a_n \,dz_0\wedge d\bar
  z_0\wedge \dotsm\wedge dz_n\wedge d\bar z_n$.  On $\Sigma_s$ we have
  $f = \bar f = 0$ and $R_s = 1$, and hence the term is equal to}
    & = -2\,\Bigl(a_0\abs{z_0}^{2(a_0-1)} + s \sum_{j=1}^n a_j
    \abs{z_j}^{2(a_j-1)} \Bigr) \, \Omega \;,
  \end{align*}
  which only vanishes, if both $s=0$ and $z_0 = 0$, i.e.\ at points
  $(0,z_1,\dotsc,z_n)\in\Sigma_0 = \widetilde\Sigma(a_0,\dotsc,a_n)$.
\end{proof}

\begin{remark}
  Note that by Lemma~\ref{isot} and~\ref{alfa_s} it follows that
  $\bigl(\Sigma(a_1,\dotsc,a_n),\alpha_1\bigr)$ is a contact
  submanifold of $\bigl(\S^{2n-1},\beta\bigr)$, and the $1$-form
  $\alpha_0$ on $\widetilde\Sigma(a_0,\dotsc,a_n)$ is equal to the
  pull-back $\pi_0^* \alpha_\mathrm{std}$ of the standard structure on
  the sphere under the projection $\pi_0:\, \C^{n+1} \to \C^n,
  (z_0,\dotsc,z_n) \mapsto (z_1,\dotsc,z_n)$.
\end{remark}

\begin{theorem}
\label{brieskorn is ctct branched cover}

  The Brieskorn manifold $(\Sigma(a_0,\dotsc,a_n), \alpha_1)$ is a contact
  branched cover of the standard sphere $(\S^{2n-1}, \alpha_\mathrm{std})$.
  More precisely: The map $\pi_0:\, \C^{n+1}\to\C^n, (z_0,\dotsc,z_n)
  \mapsto (z_1,\dotsc,z_n)$ induces an $a_0$-fold cyclic branched contact
  covering
  \begin{align*}
    \pi_0:&\, \widetilde\Sigma(a_0,\dotsc,a_n) \to (\S^{2n-1},\beta)
  \end{align*}
  with branching locus $\Sigma(a_1,\dotsc,a_n)\subset\S^{2n-1}$.
\end{theorem}
Note that the latter statement justifies the former, because
$(\S^{2n-1},\alpha_\mathrm{std}) \cong (\S^{2n-1},\beta)$, and
$(\Sigma(a_0,\dotsc,a_n),\alpha_1)$ is contactomorphic to
$\widetilde\Sigma(a_0,\dotsc,a_n)$ with the contact structure induced by
the branched covering.

\begin{proof}
  It has been known for a long time that $\pi_0$ restricted to
  $\widetilde\Sigma(a_0,\dotsc,a_n)$ is a branched covering over the
  sphere.  This can be easily seen by noting that
  $\pi_0\bigl(\widetilde\Sigma(a_0,\dotsc,a_n)\bigr)\subset\S^{2n-1}$,
  and that this map is surjective follows because a point
  $(z_1,\dotsc,z_n) \in \S^{2n-1}$ is covered by
  $(z_0,z_1,\dotsc,z_n)\in\widetilde\Sigma(a_0,\dotsc,a_1)$, where
  $z_0$ is one of the roots $\sqrt[{a_0}]{-(z_1^{a_1} + \dotsb +
    z_n^{a_n})}$.  Every point of the sphere is covered by $a_0$ points
  with the exception of the points on the branching locus
  $\Sigma(a_1,\dotsc,a_n)\subset\S^{2n-1}$.
  
  As remarked above, the $1$-form $\alpha_0$ on
  $\widetilde\Sigma(a_0,\dotsc,a_n)$ is equal to $\pi_0^*\beta$. By adding
  a small $1$-form $\epsilon\,\gamma$ to $\alpha_0$ such that
  $\left.d\gamma\right|_{\ker d\pi_0} > 0$, we obtain a contact form.  A
  possible choice for such a form is
  \begin{align*}
    \gamma &= \frac{i}{2}\,\bigl(z_0\,d\bar z_0 - \bar z_0\,dz_0\bigr)
  \end{align*}
  for sufficiently small $\epsilon > 0$, because the kernel of $d\pi_0$ is
  only non-trivial at $(0,z_1,\dotsc,z_n)
  \in\widetilde\Sigma(a_0,\dotsc,a_n)$, and the kernel lies in the
  $z_0$-plane.
  
  The only thing left to show is that
  $\bigl(\widetilde\Sigma(a_0,\dotsc,a_n),\alpha_0+\epsilon\,\gamma\bigr)$ is
  contactomorphic to $(\Sigma(a_0,\dotsc,a_n),\alpha_1)$.  This is most
  easily seen by using the contact forms $\widetilde\alpha_s =
  \left(\Phi_s^{-1}\right)^*\alpha_s$ on $\widetilde\Sigma(a_0,\dotsc,a_n)$
  for $s\in(0,1]$. This is a smooth family of forms that connects to
  $\alpha_0$, and the derivative
  \begin{align*}
    \left.\frac{d}{ds}\right|_{s=0} \widetilde\alpha_s =  \frac{i}{2}\,\bigl(z_0\,d\bar z_0 - \bar z_0\,dz_0\bigr)
  \end{align*}
  has the properties needed to apply Lemma~\ref{contact branched covering}.
\end{proof}

\begin{remark}
  It is interesting to consider, whether
  \begin{align*}
    \alpha_- := \frac{i}{2}\Bigl(-Ca_0\,(z_0\,d\bar z_0 - \bar z_0\,d z_0)
    + \sum_{j=1}^n a_j\,(z_j\,d\bar z_j - \bar z_j\,d z_j)\Bigr)
  \end{align*}
  for very large $C>0$ also gives a contact form.  The rationale is that
  the open book decomposition of such a manifold would have the same pages,
  but the monodromy map would be inverted.
  
  To check that $\alpha_-$ is a contact form, the following term should not
  vanish
  \begin{align*}
    \alpha_-\wedge \bigl(d\alpha_-)^{n-1}\wedge dR_1\wedge df \wedge d\bar
    f & = \frac{i^n(n-1)!}{2}\, a_0\dotsm a_n\,\Bigl( -2
    a_0\abs{z_0}^{2(a_0-1)} + 2C\,\sum_{j=1}^n
    a_j\abs{z_j}^{2(a_j-1)} \\
    &-(C-1)\,(a_0-1)\, \Bigl( \bar z_0^{a_0}\sum_{j=1}^n a_jz_j^{a_j}
    + z_0^{a_0}\sum_{j=1}^n a_j\bar z_j^{a_j} \Bigr)\Bigr)\,dz_0\wedge
    \dotsm\wedge d\bar z_n\;.
  \end{align*}
  It is easy to see, that this is the case for $a_0 = -1$, i.e.\ one
  gets a large set of potentially different contact structures on the
  sphere. For all Brieskorn manifolds $\Sigma(a_0,a_1,\dotsc,a_1)$, it
  is also easy to check that $\alpha_-$ is a contact form. In
  particular on $\Sigma(k,2,\dotsc,2)$, it can be shown by an explicit
  computation like the one in \cite{Koert} that the open book
  decomposition uses a $k$-fold right-handed Dehn twist for the
  monodromy map, which is indeed the inverse of the standard
  monodromy.

  Unfortunately for general combinations of integers $a_j\in\N$, it is
  quite easy to find examples where the contact condition breaks
  down.
\end{remark}

Finally, the following theorem describes a Milnor open book on $\S^{2n-1}$
which supports the contact structure $\beta$.

\begin{theorem}
  \label{milob}
  Define the polynomial $f(z_1,\dotsc,z_n) =
  z_1^{a_1}+\dotsb+z_n^{a_n}$ on $\C^n$ with $a_j\in\N$.  The sphere
  $\S^{2n-1}$ can be given an open book with binding
  $B:=\Sigma(a_1,\dotsc,a_n):=\S^{2n-1}\cap f^{-1}(0)$, and page
  fibration
  $$
  \theta:\,\S^{2n-1}-B\to\S^1,\, \z\mapsto
  \frac{f(\z)}{\abs{f(\z)}}\;,
  $$
  with $\z=(z_1,\dotsc,z_n)$.  The contact form $\beta$ on
  $\S^{2n-1}$ is supported by this open book.
\end{theorem}
\begin{proof}
  Milnor showed in \cite{Milnor_Hypersurfaces} that the structure
  defined in the lemma is an open book. Hence it only remains to show
  that $(B,\theta)$ supports the contact form $\beta$.

  The binding $B$ is a Brieskorn manifold and $\beta$ is a contact
  form for such a manifold as proved in Lemma~\ref{alfa_s}.

  To show that $d\beta$ is a symplectic form on a page
  $P_{\theta_0}=\theta^{-1}(\theta_0)$, note that the map
  $$
  e^{it}\cdot(z_1,\dotsc,z_n) = (e^{it/a_1}z_1,\dotsc,e^{it/a_n}z_n)
  $$
  is a diffeomorphism from a page $P_{\theta_0}$ to
  $P_{\theta_0+t}$, and at the same time it is the flow of the Reeb
  field $\Reeb$ of $\beta$.
  $$
  \Reeb = \frac{d}{dt} (e^{it/a_1}z_1,\dotsc,e^{it/a_n}z_n) =
  \sum_{j=1}^n \frac{1}{a_j} (x_j\frac{\partial}{\partial y_j} -
  y_j\frac{\partial}{\partial x_j})\;.
  $$
  One computes that $\iota_{\Reeb}d\beta = -2
  d\left(\sum_{j=1}^n\abs{z_j}^2\right)$, and $\beta(\Reeb) =
  \sum_{j=1}^n \abs{z_j}^2 =1$.  The Reeb field is transverse to the
  pages $P_{\theta_0}$, and hence
  $\left.d\beta\right|_{P_{\theta_0}}$ is non-degenerate.
  
  Finally if the binding $B$ is connected, the orientation of $B$ as
  boundary of the page $\overline{P}_{\theta_0}$ and as contact manifold
  $(B,\beta)$ is compatible by Remark~\ref{orientation condition for
    connected binding}.  If $B$ is non-connected (which is only the case
  for $\dim B = 1$, because $(2n+1)$-dimensional Brieskorn manifolds are
  $(n-1)$-connected) each component of $B = \Sigma(a_0,a_1)$ can be written
  in the form
  \begin{align*}
    \Bigl\{(e^{i\phi/a_0},Ae^{i\phi/a_1})\Bigm| \phi\in
    [0,2\pi\lcm(a_0,a_1)]\Bigr\}\;,
  \end{align*}
  where $A_1$ is an $a_1$-th root of $-1$. The $\phi$-parametrisation gives 
  the correct orientation, and it follows that the integral of $\alpha$
  over any of the $N$ components of $B$ has the same value $C$. In particular
  it follows,
  \begin{align*}
    0 < \int_P d\alpha = \int_B \alpha = NC\;,
  \end{align*}
  and hence $C>0$.
\end{proof}

\section{Topological description of the monodromy of the open book 
of $\Sigma(a_1,\ldots,a_n)$}

In \cite{Milnor_Hypersurfaces} Milnor worked out the topology of
the page of the above open book  $(B,\theta)$ of
$\Sigma(a_1,\dotsc,a_n)$ and described the monodromy $\psi$
by its action on $H_1(B)$. Let $\Omega_a$ denote the finite cyclic
group consisting of all $a$-th roots of unity and let $J$ denote
all linear combinations $(t_1\omega_1,\dotsc,t_n\omega_n)$ where
$\omega_i\in\Omega_{a_i}$, $t_i\geq 0$ $(1\leq i \leq n)$ and
$t_1+\dotsb+t_n=1$. Then $J$ is a deformation retract of the fiber
$P_1=\theta^{-1}(1)$ (op.cit., Lemma~9.2). Here the dimension of
$P_1$ is $2n$. Furthermore, the free Abelian group $H_{n}(P_1;\Z)$
has rank $\mu=(a_1-1)\dotsm(a_n-1)$ (op.cit., Theorem~9.1).

It is straightforward to prove the following fact which appears 
in a more general setting
in \cite[Theorem~3]{ACampo} for $n=2$.

\begin{lemma}
\label{balloon}
There is a basis for $H_{n-1}(J)$ in which the $\mu\times\mu$
matrix $\Psi$ for the monodromy $\psi$ of the open book of
$\Sigma(a_1,a_2,...,a_n)$ ($n>1$ and $\gcd(a_1,a_2,\cdots,a_n)=1$) is
$$ \Psi=A_{a_1-1}\otimes A_{a_2-1}\otimes \cdots \otimes A_{a_n-1} $$
where $A_p$ is the $p\times p$ matrix given by
$$ A_p=\left[ \ba{ccccc}
1 & 1 &1 & \ldots & 1\\
-1& 0 & 0 & \ldots & 0\\
0&-1 & 0 & \ldots & 0 \\
\vdots& & \ddots & & \vdots \\
0&0&\ldots&-1 & 0 \ea \right].
$$
\end{lemma}

As the last goal, we want to express this monodromy as a product of
Dehn twists along Lagrangian spheres. In
dimension $3$ (i.e. $n=2$), each circle is Lagrangian on the $2$
dimensional pages. Furthermore, in a rational homology sphere the
binding determines the open book decomposition up to isotopy (we
learned this from \cite{CaubelPopescu}).  Hence given the binding in
a Brieskorn sphere, any corresponding description of the monodromy in
terms of Dehn twists is the solution. This has been described in a
purely topological manner, for example, in
\cite[Theorem~1]{AkbulutOzbagci}. The question remaining is the
relation between the cycles of Dehn twists in these descriptions and 
the generators of $H_1(J)$ that appear in Lemma~\ref{balloon}.

For higher dimensions the problem is more complicated.  The skeleton given
by Milnor can be made piecewise smooth, and the smooth segments are
Lagrangian submanifolds.  Unfortunately, we do not yet know how to find
proper Lagrangian embeddings of the spheres that constitute the skeleton of
a page as a bouquet. 

\bibliography{feritklaus}

\providecommand{\bysame}{\leavevmode\hbox to3em{\hrulefill}\thinspace}
\providecommand{\MR}{\relax\ifhmode\unskip\space\fi MR }
\providecommand{\MRhref}[2]{%
  \href{http://www.ams.org/mathscinet-getitem?mr=#1}{#2}
}
\providecommand{\href}[2]{#2}
\begin{thebibliography}{CPP04}

\bibitem[A'C73]{ACampo}
N.~A'Campo, \emph{Sur la monodromie des singularit\'es isol\'ees
  d'hypersurfaces complexes}, Invent. Math. \textbf{20} (1973), 147--169.
  \MR{MR0338436 (49 \#3201)}

\bibitem[AO01]{AkbulutOzbagci}
S.~Akbulut and B.~Ozbagci, \emph{Lefschetz fibrations on compact {S}tein
  surfaces}, Geom. Topol. \textbf{5} (2001), 319--334 (electronic).
  \MR{MR1825664 (2003a:57055)}

\bibitem[Aur00]{Aur}
D.~Auroux, \emph{Symplectic $4$-manifolds as branched coverings of
  $\mathbb{C}{P}^2$}, Invent. Math. \textbf{139} (2000), no.~3, 551--602.

\bibitem[CPP04]{CaubelPopescu}
C.~Caubel and P.~Popescu-Pampu, \emph{On the contact boundaries of normal
  surface singularities}, C. R. Math. Acad. Sci. Paris \textbf{339} (2004),
  no.~1, 43--48. \MR{MR2075231 (2005e:32046)}

\bibitem[Gei97]{Geiges_Construction}
H.~Geiges, \emph{{Constructions of contact manifolds}}, Math. Proc. Camb.
  Philos. Soc. \textbf{121} (1997), no.~3, 455--464.

\bibitem[Gir02]{Giroux}
E.~Giroux, \emph{{Géométrie de contact: De la dimension trois vers les
  dimensions supérieures. (Contact geometry: From dimension three to higher
  dimensions).}}, {Li, Ta Tsien (ed.) et al., Proceedings of the international
  congress of mathematicians, ICM 2002, Beijing, China, August 20-28, 2002.
  Vol. II: Invited lectures. Beijing: Higher Education Press. 405-414 }, 2002
  (French).

\bibitem[Gon87]{Gonzalo}
J.~Gonzalo, \emph{{Branched covers and contact structures.}}, Proc. Amer. Math.
  Soc. \textbf{101} (1987), 347--352.

\bibitem[KN05]{Koert}
O.~van Koert and K.~Niederkr{\"u}ger, \emph{{Open Book decompositions for
  contact structures on Brieskorn manifolds}}, Proc. Amer. Math. Soc.
  \textbf{133} (2005), 3679--3686.

\bibitem[Mil68]{Milnor_Hypersurfaces}
J.~Milnor, \emph{Singular points of complex hypersurfaces}, Annals of
  Mathematics Studies, No. 61, Princeton University Press, Princeton, N.J.,
  1968. \MR{MR0239612 (39 \#969)}

\end{thebibliography}

\end{document}